\documentstyle [twoside]{article}
\newtheorem{thm}{Theorem}[section]
\newtheorem{prop}[thm]{Proposition}
\newtheorem{cor}[thm]{Corollary}
\newtheorem{lem}[thm]{Lemma}
\newtheorem{remark}{Remark} 
\newtheorem{dfn}{Definition} 

\def\b0{{\bf 0}}
\def\ba{{\bf a}}
\def\bi{{\bf i}}

\def\bK{{\bf K}}

\def\bz{{\bf z}}
\def\b1{{\bf 1}}
\def\bA{{\bf A}}
\def\bT{{\bf T}}
\def\bV{{\bf V}}
\def\bP{{\bf P}}

\def\bR{{\bf R}}
\def\bM{{\bf M}}
\def\bD{{\bf D}}
\def\bN{{\bf N}}

\def\bZ{{\bf Z}}
\def\bC{{\bf C}}
\def\bzeta{\zeta}

\def\Blat{\mbox{\it \raise2pt\hbox{"}\kern-2pt H}}
\def\lvup{\rlap{\ ${}^{\ell\atop{\hbox{${}^{\vee}$}}}$}\cdots}
\begin{document}
\begin{center}
{\center{\Large{\bf 
On Horn-Kapranov uniformisation of the discriminantal loci
} }}

 \vspace{1pc}
{ \center{\large{ Susumu TANAB\'E }}}
\end{center}

\noindent
\begin{center}
 \begin{minipage}[t]{10.2cm}
{\sc Abstract.} {\em
In this note we give a rational uniformisation
equation of the discriminant loci associated to a
non-degenerate affine complete 
intersection variety.
To show this formula we establish a relation of the
fibre-integral with the 
hypergeometric function of  Horn  and that of  Gel'fand-Kapranov-Zelevinski.}
 \end{minipage} \hfill
\end{center}
 \vspace{1pc}
{
\center{\section{Introduction}}
}

In this note we give a concrete rational uniformisation
equation for the discriminantal loci of non-degenerate
affine complete intersection depending 
on deformation parameters.

First of all, let us fix the situation. For the complex varieties
 $ X =$ ${\bf C^\times}^{N}$ and
$S = {\bf C}^k, $ 
we consider the mapping,
$$ f: X \rightarrow S  \leqno(0.1)$$
such that  $X_s :=\{(x_1,\cdots, x_{N})\in X ; f_1(x) +s_1=0,
\ldots, f_k(x)+ s_k=0 \}.$ Let $f_1(x), \cdots,f_k(x)$ 
be polynomials that define a non-degenerate
complete intersection (CI) in the sense of
Danilov-Khovanski \cite{DX1} with the following specific form:
$${f}_\ell(x)=
x^{\vec \alpha_{1,\ell}}+\cdots +
x^{\vec \alpha_{\tau_\ell,\ell}}, \; 1 \leq \ell \leq k, \leqno(0.2)$$
where $\vec \alpha_{i,\ell} \in (\bZ_{\geq 0})^N.$
Let $n$ be the dimension of the variety $X_0$, $dim\;X_0 
=n \geq 0.$
$W_s :=\{(x_1,\cdots, x_{N}, y_1, \cdots,$ $y_k)$ $ \in X
\times ({\bf C})^k; y_1(f_1(x) +s_1) + \ldots + y_k( f_k(x) +s_k)
=0 \}.$  Then it is known that the discriminantal loci of
$X_s$ coincides with  
that of $W_s.$ That is to say, the study of the 
discriminantal loci of a CI can be reduced to that of an hypersurface 
associated with the original CI in a special manner.
This fact has been discovered by Arthur 
Cayley \cite{GZK2} and thus the method to 
reduce the geometric study of a CI to that of a hypersurface is named 
"Cayley trick" in general, even in contexts apart from the study of discriminantal loci (e.g. the description of the mixed
Hodge structure of the former by means of the latter given by T.Terasoma, A.Mavlyutov \cite{Mav} and others). Here we return to the initial spirit of Cayley who  treated the question of the discriminantal loci.

The main idea is based on that of the paper \cite{Horn} which states that
the singular loci of the linear differential operators 
annihilating the fibre integrals
of $X_s$ coincide with the discriminantal loci of $X_s.$
In the modern terminology of the  A-hypergeometric functions (HGF), 
it is equivalent to say that
A-discriminantal loci are singular loci for generalized A-HGF. This fact 
has been proven in \cite{Kap1} and we give a more precise description 
of the discriminantal loci by means of combinatorial data of the polynomial mapping $f$ and the toric geometry of $W_s$ (see Theorem  2.6).

%
 \vspace{1.5pc}
\footnoterule

\footnotesize{AMS Subject Classification: 14M10, 55R80, 32S40, 33C65.

Key words and phrases: complete intersection, discriminantal set, Horn hypergeometric function.}
\normalsize

\newpage

\par
Let us review the contents of the note in short.
In \S 1 we recall some basic facts on the Cayley trick and
N\'eron-Severi torus. 
In \S 2, we calculate the Mellin transform of the fibre integral
in an explicit manner. Using a representation of 
the Mellin transform 
we show that fibre integral satisfies the Horn type system of differential equations (Theorem 2.4). From this expression of the Horn type system, we get the discriminantal loci as the boundary of a convergence domain of solutions
 to the system.
In \S 3, we show that the fibre integral calculated in
\S 2 is nothing but the quotient of the Gel'fand-
Kapranov-Zelevinski generalized
hypergeometric function (HGF) by the torus 
action. In \S 4 we give two computational examples: 
discriminantal loci for the $D_4$ type singularity
and the simplest non-quasihomogeneous complete intersection.

Finally we remark that this note is an abridged version of some parts from 
$\cite{Tan}$ where one can find more details. 

 \vspace{2pc}
{
\center{\section{ Cayley trick and N\'eron-Severi torus}}
}

Throughout this section we keep the notation of \S 0.
Further we introduce the following notations.
Let  ${\bf T}^m=( {\bf C}\setminus \{0\})^m$ $=(\bC^{\times})^m$
be the complex algebraic torus of dimension $m.$ We denote by
$x^{\bi}$ the monomial $x^{\bi}:= x_1^{i_1} \cdots x_N^{i_N}$ 
with multi-index ${\bi}=({i_1}, \cdots, {i_N}) \in {\bf Z}^{N},$ and by
$dx$ the $N-$volume form $dx :=dx_1 \wedge\cdots \wedge dx_N.$ 
We shall also use the notations
$x^\b1 :=x_1\cdots x_N, $ $ y^{\bzeta}= y_1^{\zeta_1}\cdots
y_k^{\zeta_k},$ $s^\bz = s_1^{z_1} \cdots s_k^{z_k}$ and $ds =
ds_1\wedge \cdots \wedge ds_k$ and their analogies for each variable. 
In this section we consider an extension of
 the mapping $f$ to that defined from
$\bP_{\tilde \Sigma}$ to $\bC^k.$ 
We follow the construction by
\cite{BC} and  \cite{Mav}.  Let us define $M$ as the dimension 
of a minimal ambient space so that we can quasihomogenize simultaneously
the polynomials
$(f_1(x),\cdots,f_{k}(x))$ by multiplying certain terms by new variables:  
$$ x^\bi
\longmapsto x'_j x^\bi, \; j=1,2, \cdots.  $$ Let us denote by
$(f_1(x,x'),\cdots,f_{k}(x,x'))$ the new 
polynomials obtained in such a way.  
These polynomials are quasi-homogeneous with respect to certain weight
system i.e. there exists a set of positive integers
$(w_1, \cdots, w_N, w_1', \cdots, w'_{M-N})$ such that their
G.C.D. equals 1 and the following relation holds:
$$ E(x,x')(f_\ell(x,x')) = p_\ell f_\ell(x,x')
\;\;\;\mbox{for}\;\ell=1,\cdots, k,$$
where $p_\ell$ is some positive integer and
$$ E(x,x')= \sum_{i=1}^N w_ix_i \frac{\partial}{\partial x_i}
+ \sum_{j=1}^{M-N}w'_j x'_j \frac{\partial}{\partial x'_j}
,\leqno(1.1)$$
$E$  an Euler vector field.

\par {\bf Example} We modify the polynomial
$ f(x)=x_1^a + x_1x_2 + x_2^b,$
with $a,b >2,$ $GCD(a,b)=1,$ in adding a new variable $x'_1$ so that the new polynomial
$ f(x,x')=x_1^a + x'_1x_1x_2 + x_2^b, $
becomes quasihomogeneous with respect to the weight system 
$(b,a, ab-a-b)$.

In general there are of course many choices of terms that we modify
to realize the quasihomogeneiety.
\par From now on we will use the notation
$X:= (X_1 , \cdots, X_M):= (x_1, \cdots, x_N,$$
x_1', \cdots,$
$ x'_{M-N})$  and that of the polynomial $f_\ell(X):= f_\ell(x,x').$
If we introduce the Euler vector field,
$$ E(X')= \sum_{i=1}^Nw_ix_i \frac{\partial}{\partial x_i}
+ \sum_{j=1}^{M-N}w'_j x'_j \frac{\partial}{\partial x'_j}+
X_{M+1} \frac{\partial}{\partial X_{M+1}},$$
we have the following relation:
$$ E(X')(f_\ell(X)+X_{M+1}^{p_\ell}s_\ell) =
p_\ell (f_\ell(X)+X_{M+1}^{p_\ell}s_\ell)
\;\;\;\mbox{for}\;\ell=1,\cdots, k.$$ From now on we denote 
$X':=(X,X_{M+1}).$
 Let $\bM_\bZ$ be an integer 
lattice of rang $N$ and  $\bN_\bZ$ be its dual, 
$\bN_\bZ= Hom(\bM_\bZ, \bZ).$
We denote by $\bM_\bR$ (resp. $\bN_\bR$) the natural extension
of $\bM_\bZ$ (resp. $\bN_\bZ$) to its real space.
Let us take
$\vec{e}_1, \cdots, \vec{e}_{M+1}$ a set of generators of one dimensional
cones such that $ \sum_{\ell=1}^{M+1}\bR \vec{e}_\ell = \bN_\bR.$
We can define a simplicial fan $\Sigma$
in $\bN_\bR$ as a set of simplicial cones spanned by
the above $\vec{e}_1, \cdots, \vec{e}_{M+1}$.
Our construction of the Euler vector field $E(X')$
correspond to the superstructure  $\bN_\bR \times \bN_\bR'$
with a basis of generators 
$\vec{\tilde e}_{N+1}, \cdots \vec{\tilde e}_{M+1}$ such that
$$ \sum_{i=1}^Nw_i\vec{\tilde e}_i +  \sum_{j=1}^{M-N}w'_{j} \vec{\tilde
e}_j +\vec{\tilde e}_{M+1}=0.$$ Here we have
$ p_\bN(\vec{\tilde e}_j) =\vec{e}_j $
for the projection
$p_\bN: \bN_\bR \times \bN_\bR' \rightarrow \bN_\bR.$
While the dimension of the vector space
$\bN_\bR \times \bN_\bR'$ must be minimal i.e.
$dim(\bN_\bR \times \bN_\bR')=$
$M.$

We introduce a polynomial,
$$H(x,y):=y_1f_1(x)+\cdots + y_k f_k(x) \in \bZ[x_1, \cdots, x_N,
y_1, \cdots, y_k], \leqno(1.2) $$
in adding new variables $y_1, \cdots, y_k.$
Let $ \vec{n}_1, \cdots, \vec{n}_{M+k}$
be the elements of the set $ supp( H(x,y)) \subset \bZ^{N+k}.$
We define a simplicial rational fan $\tilde \Sigma$
in $\bR^{N+k}$ as a set of simplicial cones generated by
 $\vec{n}_1, \cdots, \vec{n}_{M+k}.$
We consider the injective homomorphism
$$ \varphi: {\tilde \bM}_\bZ \rightarrow \bZ^{M+k},$$
for ${\tilde \bM}_\bZ = \bM_\bZ \times \bZ^k,$ defined by
$$\varphi(\vec{\tilde m}) =(<\vec{\tilde m}, {\vec{n}_1}>, \cdots,
<\vec{\tilde m}, {\vec{n}_{M+k}}>).$$
The cokernel of this mapping is a free abelian group,
$$ {  C l}( \tilde \Sigma) = \bZ^{M+k}/\varphi(\tilde \bM_\bZ)$$
for which the following group can be defined
$$\bD (\tilde \Sigma) := Spec \bC [{  C l}(\tilde \Sigma) ]. \leqno(1.3)$$
As a matter of fact  this group $\bD (\tilde \Sigma)$ is isomorphic to 
an algebraic torus $\bT^{M-N}.$
One can define the toric variety $\bP_{\tilde \Sigma}$
associated to the affine space,
$$ \bA^{M+k}= Spec \bC [X_1, \cdots, X_M,
y_1, \cdots, y_k ].$$
To this end we proceed following way after the method initiated by
M.Audin.
Let  $\hat{X_\sigma} := \prod_{1 \leq i \leq 
M, \vec{n}_i \not \in  \sigma} X_i$
$\prod_{1 \leq j \leq k, \vec{n}_{M+j} \not \in  \sigma} y_j,$
be a monomial defining a coordinate plane and the ideal
$$B(\tilde \Sigma) = <\hat{X_\sigma}; \sigma \in \tilde \Sigma> 
\subset \bC [X_1, \cdots, X_{M}, y_1, \cdots, y_k].$$
Let $ Z(\tilde \Sigma):=\bV(B(\tilde \Sigma)) \subset \bA^{M+k} $
be the variety defined by the ideal $B(\tilde \Sigma).$
We construct the toric variety $\bP_{\tilde \Sigma}$ 
as the quotient of $U(\tilde \Sigma):=\bA^{M+k} \setminus
Z(\tilde \Sigma)$ by the group action $\bD(\tilde \Sigma)$:
$$ \bP_{\tilde \Sigma} = U(\tilde \Sigma )/ \bD(\tilde \Sigma),$$
with $ dim \;\bD(\tilde \Sigma)= M-N, dim\;U(\tilde \Sigma )= M+k.$
\begin{dfn}
This group $\bD(\tilde \Sigma) \cong \bT^{M-N}$ is called
the N\'eron-Severi torus associated to the fan $\tilde \Sigma.$
\label{dfn1}
\end{dfn}

\par

We introduce the following polynomial (named phase function below),
$$F(X,s,y):=y_1(f_1(X)+s_1)+\cdots + y_k (f_k(X)+s_k), \leqno(1.4)$$
that will play essential r\^ole in our further studies.
In \S 3, we treat the following affine variety defined for $(1.4)$:
 $$ Z_{F(x,\b1,\b1,y)+1} =\{(x,y) \in \bT^{N+k}; F(x,\b1,\b1,y)+1=0 \}. 
\leqno(1.5)$$

Further on we shall prepare several lemmata on combinatorics
which are useful for the derivation of the discriminant loci 
equation.
We denote by $L$  the number of monomials in
$\;(X,s,y)\;$ that take part in the phase function $(1.4)$ for $(0.2)$. 
That is to say $L= \sum_{q=1}^k (\tau_q+1).$
Here we introduce new variables $(T_1, \cdots, T_L)$ $\in$
$\bT^L$ that satisfy the following relations,
$$ T_1 = y_1 x^{\vec \alpha_{1,1}},T_2 = y_1 x^{\vec \alpha_{2,1}}, \cdots,
T_{L} = y_k s_k. \leqno(1.6)$$ 
Each $T_q$ represents the $q-$th monomial present in $F(x,\b1,s,y)$
(see (2.3) below).
 We will use the following matrix ${\sf M}(A)$
whose column is a vertex of the Newton polyhedron
$\Delta(F(x,\b1,\b1,y)),$ 
$${\sf M}(A):=
 \left [\begin {array}{cccccccccccccc}
1&1&\cdots&1&0&0&\cdots&0&0&\cdots&0&0&\cdots&0\\
0&0&\cdots&0&1&1&\cdots&1&0&\cdots&0&0&\cdots&0\\
0&0&\cdots&0&0&0&\cdots&0&1&\cdots&0&0&\cdots&0\\
\vdots&\vdots&\vdots&\vdots&\vdots&\vdots
&\vdots&\vdots&\vdots&\vdots&\vdots&\vdots&\vdots&\vdots\\
0&0&\cdots&0&0&0& \cdots&0&0&\cdots&1&1&\cdots&1\\
0&\alpha_{111}&\cdots&\alpha_{\tau_111}&0&\alpha_{121}& \cdots&
\alpha_{\tau_221}&0&\cdots&0&\alpha_{1k1}&\cdots&\alpha_{\tau_kk1}\\
\vdots&\vdots&\vdots&\vdots&\vdots&\vdots
&\vdots&\vdots&\vdots&\vdots&\vdots&\vdots&\vdots&\vdots\\
0&\alpha_{11N}&\cdots&\alpha_{\tau_11N}&0&\alpha_{12N}&\cdots&
\alpha_{\tau_22N} &0&
\cdots&0&\alpha_{1kN} &\cdots&\alpha_{\tau_kkN}\\ \end {array}\right ].
\leqno(1.7) $$
Further we assume that $rank({\sf M}(A))= k+N.$ 
We always assume the inequality $N+2k \leq L$ for $(0.2)$.

In this situation we can define a non-negative integer
$m$ as the minimal number of variables
$$x'' =(x'_1, \cdots, x'_{m}) \leqno(1.8)$$
to make the number of variables present in the expression $(1.4)$
equal to $L$. That is to say $L= N+m+2k.$
For example, the relation $(1.6)$ may be modified into the following form:
$$ T_1 = y_1x_1' x^{\vec \alpha_{1,1}},
T_2 = y_1x_2' x^{\vec \alpha_{2,1}}, \cdots,
T_{L-1} = y_k x'_m x^{\vec \alpha_{\tau_k,k}},
T_{L} = y_k s_k. \leqno(1.6)'$$
In other words, proper addition 
of new variables $x''=
(x'_1, \cdots, x'_{m})$
to  $f_1(x), \cdots, f_k(x)$  makes the polynomial $F(X,0,y)$
quasihomogeneous. In this way we have
$$M= N+m.   \leqno(1.9)$$

Further we shall consider a simple parametrisation of the variety
$$ Z_{F(X,s,y)} =\{(X,y) \in \bT^{M+k}; F(X,s,y)=0 \}. 
\leqno(1.10)$$
 
Namely we denote,
$$ \Xi := ^t(x_1, \cdots, x_N,  x_1', \cdots,  x'_m,  s_1, \cdots,
 s_k,y_1, \cdots,   y_k) , \leqno(1.11)$$
$$Log\; T := ^t(log\; T_1, \cdots, log\; T_L)
\leqno(1.12)$$
$$Log\; \Xi 
:= ^t(\log\; x_1, \cdots, \log\; x_N,   log\; x_1', \cdots, log\; x'_m,
\log\; s_1, \cdots,
log\; s_k,\log\; y_1, \cdots, log\; y_k). \leqno(1.13)$$ 
Then we have, for example, a linear equation equivalent to $(1.6)'$
that can be written down as follows,
$$ log\; T_1 = log\; y_1 + log\; x_1' +<\vec {\alpha_{1,1}}, log\; x>,
log\; T_2= log\; y_1+ log\; x_2' + <\vec {\alpha_{2,1}} log\; x>, \cdots,
\leqno(1.14)$$
$$
log\; T_{L-1} = log\; y_k + log\; x'_m+ <{\vec \alpha_{\tau_k,k}}, log\; x>,
log\; T_{L} = log\; y_k + log\; s_k. $$

Let us write down the relation between $(1.12)$ and $(1.13)$ 
by means of a matrix
${\sf L} \in End({\bf Z}^L),$
$$ Log\; T= {\sf L}\cdot Log\; X. \leqno(1.15)$$
Below the columns $\vec v_i$ (resp. $\vec w_i$)
of the matrix ${\sf L}$ (resp. ${\sf L}^{-1} $) 
shall always be ordered in accordance with
(1.11), (1.12), (1.13) unless otherwise stated.

For the polynomial mapping (0.2), the choice of monomials to be modified by
supplementary variables is a bit delicate.
Namely, we have to observe the following rules to avoid the
degeneracy of the matrix  $\sf L$ of the relation
$(1.15).$
\begin{lem}
For (0.2) and $(1.8)$, we get a non-degenerate matrix $\sf L$
if we observe the following rules:
\par {\bf a.} For the fixed index $q \in\{1,\cdots, k \}$, 
it is necessary to choose at least one of monomials 
$x^{\vec\alpha_{i,q}}, 1 \leq i \leq \tau_q$ that remains without modification.
\par {\bf b.} For the fixed index $j \in\{1,\cdots, N\}$
it is necessary to choose at least one of monomials 
$x^{\vec\alpha_{r,i}}$ such that $\alpha_{r,i,j}\not =0$,  
$ 1 \leq i \leq k,$ $ 1 \leq r \leq \tau_i,$  
that remains without modification.
\label{lem121}
\end{lem}

We recall here the notion of non-degenerate hypersurface,
\begin{dfn}
The hypersurface defined by a polynomial $g(x)= \sum_{\alpha \in supp(g)}
g_\alpha x^\alpha$ $\in \bC[x_1, \cdots,
x_n]$ is said to be non-degenerate if and only if for any
$\xi \in \bR^n$ the following inclusion takes place,
$$ \{x \in\bC^n;
x_1 \frac{\partial
g^\xi}{\partial x_1}= \cdots=  x_n \frac{\partial g^\xi}{\partial
x_n}=0 \}\subset \{x \in\bC^n;x_1 \cdots x_n =0\}$$
where $g^\xi(x)= \sum_{\{\beta; <\beta, \xi> \leq <\alpha, \xi>,
\;{\rm for \;\;all}\;\;\alpha \in supp(g) \}} g_\alpha x^\alpha.$
We call the CI $X_0$ for (0.2) non-degenerate if the hypersurface
$Z_{F(x,\b1,0,y)+1}$ is non-degenerate.
\label{dfn15}
\end{dfn}
The following is an easy consequence of the above Definition.
\begin{prop}
If the matrix $\sf L$ is non-degenerate, the hypersurface
$Z_{F(x,\b1,0,y)+1}$ and the CI $X_0$ are 
non-degenerate in the sense of the Definition
~\ref{dfn15}.
\label{prop18}
\end{prop}

{
\center{\section{ Horn's hypergeometric functions }}
}

From this section, we change the name of variables $x"=(x'_1, \cdots, x'_m)$
into $s':=(s'_1, \cdots, s'_m).$ We use both of the notations 
$X=(x,x")= (x,s').$

Let us consider the Leray's coboundary (see \cite{Vas}) to define the
fibre integral,  
$\gamma \subset H_{N}({\bf T}^{N} \setminus
\cup_{i=1}^k \{x \in \bT^N: f_i(X)+s_i=0\})$ such that
$\Re (f_i(X) + s_i)|_{\gamma} <0$.
Further on central object of our study is the following fibre integral,
$$  I_{x^{\bi}, \gamma}^{\zeta}(s,s')=
\int_{\gamma}
 (f_1(x,{  s}') +s_1)^{-\zeta_1-1}
\cdots (f_k(x,{s}') +s_k)^{-\zeta_k-1} x^{\bi+\b1} \frac{dx}{x^{\b1}},
\leqno(2.1)$$ and its Mellin transform,
$$ M_{{\bi},\gamma}^\zeta ({\bz},{\bz'} ):=\int_\Pi s^{\bz }s'^{\bz'}
I_{x^{\bi}, \gamma}^{\zeta} (s,s')
\frac{ds}{s^{\b1}}\wedge \frac{ds'}{s'^{\b1}},
\leqno(2.2)$$ 
for  certain cycle $\Pi$ homologous to ${\bf R}^{m+k}$
which avoids the singular loci of 
$I_{x^{\bi}, \gamma}^{\bzeta} (s,s')$  (cf. \cite{Tsikh}).
After Definition ~\ref{dfn1} above, we understand that
$s' \in \bD(\tilde \Sigma)$ is a variable on the 
N\'eron-Severi torus. Thus the fibre integral 
$I_{x^{\bi}, \gamma}^{\zeta}(s,s')$ is a ramified function on the torus
$\bD(\tilde \Sigma) \times \bT^k.$
It is useful to understand the calculus of the Mellin
transform 
in connection with the notion of the generalized
HGF in the sense of Mellin-Barnes-Pincherle
\cite{AK}, \cite{Nor}.
After this formulation, the classical HGF of Gauss can be expressed 
by means of the integral,
$$ _2 F_1(\alpha,\beta,\gamma|s)=
\frac {1}{2\pi i } \int_{z_0 - i\infty}^{  z_0 + i\infty }(-s)^z
\frac {\Gamma ( z+ \alpha)\Gamma ( z+ \beta )\Gamma ( -z)}{\Gamma
( z+ \gamma)} dz , \;\; - \Re  \alpha, - \Re  \beta < z_0.$$

Next we modify the Mellin transform
$$
M_{{\bi},\gamma}^\zeta ({\bz}, \bz' )
= c(\zeta) \int_{S^{k-1}_+(w'') \times \gamma^\Pi} \frac{x^{\bi}
\omega^\zeta s^{\bz -\b1} s'^{\bz' -\b1} dx \wedge \Omega_0(\omega) 
\wedge 
ds \wedge ds'}
{(\omega_1(f_1(X) +s_1)+\cdots +\omega_k(f_k(X) +s_k))^ {\zeta_1+
\cdots +\zeta_k +k}} $$ $$= c(\zeta)\int_{\bR_+} \sigma^{\zeta_1+ \cdots
+\zeta_k +k} \frac{d \sigma}{\sigma} \int_{S^{k-1}_+(w'')}
\omega^\zeta \Omega_0(\omega) \int_{\gamma}x^{\bi} dx \int_{\Pi}
s^{\bz}s'^{\bz'} e^{\sigma(\omega_1(f_1(X) +s_1)+\cdots 
+\omega_k(f_k(X)
+s_k))} \frac{ds}{s^\b1}\frac{ds'}{{s'}^\b1},$$
with $ c(\zeta) = \frac{\Gamma(\zeta_1+\cdots +\zeta_k+k)}
{\Gamma(\zeta_1+1)\cdots \Gamma(\zeta_k+1)}.$ Here we made 
use of notations
$S^{k-1}_+(w'') =\{(\omega_1, \cdots,
\omega_k): \omega_1^{\frac{\bf w''}{w_1''}}+\cdots
+\omega_k^{\frac{\bf w''}{w_k''}}=1, \omega_\ell >0\;\;$ 
for all $\; \ell,$ ${\bf w''}$ $=$ $\prod_{1 \leq i \leq k}w''_i \}$ and
$\Omega_0(\omega)$ the $(k-1)$ volume form on $S^{k-1}_+
(w''),$
$$\Omega_0(\omega) = \sum_{\ell=1}^k (-1)^\ell w''_\ell \omega_\ell d\omega_1\wedge \lvup \wedge
d\omega_k.$$
In the above transformation we used a classical interpretation of
Dirac's delta function as a residue:
$$\int_{\gamma} \int_{\bR_+}
e^{ y_j(f_j(X) + s_j)} y_j^{\zeta_j} dy_j \wedge dx
= \Gamma(\zeta_j+1) \int_{\gamma}
 (f_j(X) + s_j)^{-\zeta_j-1} dx. $$

We introduce the  notation $\gamma^\Pi:= \cup_{(s,s')\in
\Pi}((s,s'), \gamma).$ One shall not confuse it with the 
thimble of Lefschetz, because  $\gamma^\Pi$ is rather a tube without thimble.
We will rewrite  the last expression,
$$\int_{(\bR_+)^{k} \times \gamma^\Pi}
e^{\Psi(T)} x^{\bi+\b1} y^{\zeta+\b1}s^{\bz} {s'}^{\bz'} \frac{dx}{x^\b1}
\wedge \frac{dy}{y^\b1} \wedge \frac{ds}{s^\b1}\wedge \frac{ds'}{{s'}^\b1}    
$$
where
$$\Psi(T) = T_1(X,s,y) + \cdots + T_L(X,s,y)= F(X,s,y),  \leqno(2.3)$$
in which each term $T_i(X,s,y)$ stands for a monomial in 
variables
$(X,s,y)$ of the phase function $(1.4).$
We transform the above integral into the following form,
$$
\int_{({\bR_+})^k \times  \gamma^\Pi } e^{\Psi(T(X,s,y))}
x^{\bi+\b1 } s^\bz s'^{\bz'} y^{{\bzeta}+\b1} \frac{dx}{x^\b1} \wedge
\frac{dy}{y^\b1} \wedge 
\frac{ds}{s^\b1}\wedge \frac{ds'}{{s'}^\b1} \leqno(2.4)$$
$$
= (det {\sf L})^{-1}\int_{{\sl L}_\ast ({ \bR_+}^k
\times \gamma^\Pi  )}
e^{\sum_{a \in I}T_a} \prod_{a\in I}
T_a^{{\mathcal L}_a({\bi, \bz, \bz', \bzeta})}
\bigwedge_{a \in I} \frac{dT_a}{T_a}  $$
$$
= (-1)^{\zeta_1+\cdots +\zeta_k +k}
(det {\sf L})^{-1}\int_{-{\sl L}_\ast ({ \bR_+}^k
\times \gamma^\Pi  )}
e^{-\sum_{a \in I}T_a} \prod_{a\in I}
T_a^{{\mathcal L}_a({\bi, \bz, \bz', \bzeta})}
\bigwedge_{a \in I} \frac{dT_a}{T_a}.  $$
Here ${\sl L}_\ast ({\bR_+}^k
\times \gamma^\Pi)$ means the image of the chain in
$\bC^{M}_X\times \bC^{k}_s\times \bC^{k}_y$
into that in  $\bC^{L}_T$
induced by the  transformation $(1.15).$
We define $-{\sl L}_\ast ({ \bR_+}^k
\times \gamma^\Pi  )$ $=\{(-T_1, \cdots, -T_L)$ $\in \bC^L;$ 
$(T_1, \cdots,T_L)$ $\in
{\sl L}_\ast ({ \bR_+}^k
\times \gamma^\Pi  ), \Re T_a <0, a \in [1,L] \}.$
The second equality of $(2.4)$ follows from Proposition 2.1, 3) below
that can be proven in a way independent of the argument to derive $(2.4).$
We will denote the set of columns and rows of the matrix 
$\sf L$ by  $I,$
$$ I:= \{1, \cdots, L\}.$$ Here we remember the  relation
$L=N+m+2k= M+2k.$

The following notion helps us to formulate the result in a compact manner.
\begin{dfn}
A meromorphic function $g(\bz,\bz')$
is called $\Delta-$periodic for  
$\Delta \in \bZ_{>0},$ if
$$g(\bz,\bz')= h(e^{2 \pi \sqrt -1
\frac{z_1}{\Delta}}, \cdots,
e^{2 \pi \sqrt -1
\frac{z_k}{\Delta}}, 
e^{2 \pi \sqrt -1
\frac{z_1'}{\Delta}}, \cdots,
e^{2 \pi \sqrt -1
\frac{z_m'}{\Delta}}),$$
for some rational function $h(\zeta_1, \cdots, \zeta_{k+m}).$
\end{dfn}

For the simplicial CI  $(0.2)$ (i.e. we can construct $F(X,s,y)$ for which
the matrix $\sf L$ is non-degenerate), we have the following statement.
\begin{prop}
1)For any cycle $\Pi \in H_{k+m}({\bf T}^{k+m} \setminus S.S.
I_{x^{\bi}, \gamma}^{\bzeta} (s,s'))$ 
the Mellin transform $(2.1)$
can be represented as a product of $\Gamma-$ function factors
up to a $\Delta-$periodic function factor $g(\bz)$,
$$  M_{{\bi}, \gamma}^\zeta (\bz,\bz' )= g(\bz) \prod_{a \in I}
\Gamma\bigl( {\mathcal L}_a({\bi, \bz,\bz', \bzeta})\bigr),$$ with
$${\mathcal L}_a({\bi, \bz, \bz',\bzeta} ) =\frac{\sum_{j=1}^N A_j^a (i_j+1)
+\sum_{j=1}^m C_j^a z'_j
+\sum_{\ell=1}^k \left( B_\ell^a z_\ell + D_\ell^a(\zeta_\ell+1)\right)
}{\Delta}, a \in I. \leqno(2.5)$$
Here the following matrix $\Delta^{-1}{\sf T} = (\sf L)^{-1}$ 
has integer elements,
$$^t{\sf T}=(A_1^a,\cdots,
A_{N}^a, C_1^a , \cdots, C_{m}^a, B_1^a
, \cdots, B_k^a ,D_1^a
, \cdots, D_k^a )_{1 \leq
a \leq L}, \leqno(2.6)$$
with $G.C.D.(A_1^a,\cdots,
A_{N}^a, C_1^a , \cdots, C_{m}^a, B_1^a
, \cdots, B_k^a ,D_1^a
, \cdots, D_k^a )=1,$
for all $ a \in [1, L].$
In this way $\Delta >0$ is uniquely determined.
The coefficients of  (2.5) satisfy the following properties for each index
 $a \in I$ :

$\bf a$
Either
${\mathcal L}_a({\bi, \bz, \bz', \bzeta} ) = \frac{\Delta}{\Delta}z_\ell,$
i.e. $A_1^a=\cdots=
A_{N}^a=0,$ $B_1^a = \lvup = B_k^a=0,$ $B_\ell^a=1.$

$\bf b$
Or
$${\mathcal L}_a({\bi, \bz, \bz', \bzeta} )=
\frac{\sum_{j=1}^N A_j^a (i_j+1)+\sum_{j=1}^{m} C_j^a z'_j
+\sum_{\ell=1}^k B_\ell^a (z_\ell -\zeta_\ell-1) }{\Delta}$$

2) For each fixed index $1 \leq \ell \leq N, 1 \leq q \leq k,$
$1 \leq j \leq m$ the following equalities take place:
$$\sum_{a \in I} A_\ell^a =0,\; \sum_{a \in I}
B_q^a =0,\sum_{a \in I}
C_j^a =0.\leqno(2.7)$$

3) The following relation holds among the linear functions
 ${\mathcal L}_a,$ $a \in I$:
$$ \sum_{a \in I }{\mathcal L}_a(\bi, \bz, \bz', \zeta)= \zeta_1 +\cdots +
\zeta_k +k.$$
\label{prop21}
\end{prop}

{\bf Proof}

1) First of all we recall the definition of the
$\Gamma-$function,
$$ \int_{C_a} e^{-T_a}T_a^{\sigma_a} \frac{dT_a}{T_a} = (1-
e^{2 \pi i \sigma_a})
\Gamma(\sigma_a),$$
for the unique non-trivial cycle $C_a$ that turns around
$T_a=0$ with the  asymptotes
$\Re T_a \rightarrow + \infty.$
We consider a transformation of  the integral $(2.4)$ 
induced by the change of cycle $\lambda:C_a \rightarrow \lambda(C_a)$
defined by the relation,
$$ \int_{\lambda(C_a)} e^{-T_a}T_a^{\sigma_a} \frac{dT_a}{T_a}
= \int_{C_a} e^{-T_a}(e^{2\pi  \sqrt -1 }T_a)^{\sigma_a} \frac{dT_a}{T_a}.$$
By the aid of this action the chain ${\sf L}_\ast (
{\bR_+}^k \times
\gamma^\Pi ) $ turns out to be homologous to a chain,
$$\sum_{(j_1^{(\rho)},\cdots,j_L^{(\rho)} )\in [1,\Delta]^L}
m_{j_1^{(\rho)},\cdots,j_L^{(\rho)}} \prod_{a=1}^k
\lambda^{j_a^{(\rho)}}(\bR_+)\prod_{a'=k+1}^L
\lambda^{j_{a'}^{(\rho)}}(C_{a'}),$$ with $m_{j_1^{(\rho)},
\cdots, j_L^{(\rho)}}  \in \bZ. $ 
This fact explains the appearance of the factor
$g(\bz, \bz')=$ $\sum_{(j_1^{(\rho)},\cdots,j_L^{(\rho)} )\in
[1,\Delta]^L}$ $m_{j_1^{(\rho)},\cdots,j_L^{(\rho)}}$ $
\prod_{a=1}^k$ $ e^{2 \pi \sqrt -1 j_a^{(\rho)}{\mathcal
L}_a({\bi, \bz, \bz', \zeta})}$ $ \prod_{a'=k+1}^L$ $ e^{2 \pi \sqrt -1
j_{a'}^{(\rho)}{\mathcal L}_{a'}({\bi, \bz, \bz', \zeta})} (1- e^{2 \pi
\sqrt -1 {\mathcal L}_{a'}({\bi, \bz, \bz', \zeta} ) })$ apart from the
factors of type  $\Gamma (\bullet).$

In the sequel we analyze the $\Gamma-$ function factors that arise 
from the integral
(2.4).
To this end, we represent the matrix ${\sf L}$ (resp.${\sf L}^{-1}$)
as a set of $L$ columns properly ordered:
$${\sf L} =({\vec v}_1,{\vec v}_2,, \cdots,{\vec v}_L),
{\sf L}^{-1} =({\vec w}_1,{\vec w}_2,, \cdots,{\vec w}_L), {\vec
w}_a =^t({w}_{a,1}, \cdots, {w}_{a,L}). \leqno(2.8)$$ 

The interior product of vectors $(\bi+\b1, \bz, \bz', \bzeta+\b1)$ 
and ${\vec w}_a$ defines the linear function in question:
$${\mathcal L}_a({\bi, \bz, \bz', \bzeta})= (\bi +\b1, \bz, \bz', \bzeta +\b1)
\cdot {\vec w}_a. \leqno(2.9)$$ 

The vector columns of 
${\sf L}^{-1} $ are divided into 3 groups:
\par {\bf 1 } the columns with
all formally non-zero elements.
\par {\bf 2} with unique non-zero element $(=1)$
that produces  $z_i, 1 \leq i \leq k$ and $z'_j, 1 \leq j \leq m$  in (2.9).
\par {\bf 3} with the non-zero elements that produce a function
linear in $\zeta +\b1, \bi +\b1$
after (2.5).
\par
In the further argument, only the first two groups of columns are  important.

The column that corresponds to $log\;s_i$ of $\sf L$ 
contains the unique non-zero element $(=1)$ at the position 
$\tau_1+ \cdots +\tau_i+i.$
Meanwhile the column of $\sf L$ that corresponds to the variable $log
\;x'_\ell$ consists also of an unique non-zero element $(=1)$ outside 
the positions $\tau_1+ \cdots +\tau_i+i, (1 \leq i \leq
k).$
Let us denote this correspondence by
$$ {\vec v}_{\rho(i)}= ^t(0,\cdots, 0,
{\rlap{\ ${}^{\sigma(i)\atop{\hbox{${}^{\vee}$}}}$}\cdots} 1, 0,
\cdots,0),$$ that yields in  $\sf L^{-1},$
$$ {\vec w}_{\sigma(i)}= ^t(0,\cdots, 0,
{\rlap{\ ${}^{\rho(i)\atop{\hbox{${}^{\vee}$}}}$}\cdots} 1, 0,
\cdots,0).$$ Here the mappings  $\rho,\sigma: \{N+1, \cdots , M+k\} \rightarrow
I$ are injections that send the number of columns corresponding 
to the variables $s,x'$ to the total set of 
indices $I$. We divide the columns of  $\sf L^{-1}$ into
$k$ groups $\Lambda_1, \cdots,\Lambda_k \subset I$ each of which 
corresponds to  $\Lambda_b = \{\tau_1+\cdots +\tau_{b-1}+b ,
\cdots,\tau_1+\cdots + \tau_{b}+b\} \subset I.$ For this group, one can 
claim following assertions. $a)$ The column
${\vec v}_{M+k+b}= ^t(0,\cdots, 0,0, {\rlap{\ ${}^{\tau_1+\cdots
+\tau_{b-1}+b \atop{\hbox{${}^{\vee}$}}}$}\cdots}
0,1,1,\cdots,\cdots, {\rlap{\ ${}^{\tau_1+\cdots + \tau_{b}+b
\atop{\hbox{${}^{\vee}$}}}$}1}, \cdots,1, 0, \cdots,0),$ with
$\tau_{b}+1,$ $(1 \leq b \leq k)$ non-zero elements $(=1).$
$b)$ For the vectors  $\vec w_a$ of the case $\bf 1$ above,
$$\sum_{a \in \Lambda_b}w_{a,j}= 0\;\;{\rm if}\; j \not= M+k+b, 1 \leq b \leq k, \leqno(2.10)$$
and there exists another vector of the same group $\Lambda_b$
that satisfies:
$$w_{\sigma(i),j}=\delta_{\rho(i),j},\;\;\leqno(2.11)$$
where $\delta_{\cdot, \ast}$ is the Kronecker delta symbol.
The vector $(2.11)$ corresponds to the group {\bf 2}.

Thus the columns of the group $\bf 2$ (resp. $\b1$) give rise to the linear functions of the group $\bf b$ (resp. $\ba$).

2) The  1-st $, \cdots,$ ${M+k}-$th vector rows of the matrix
$\sf L^{-1}$ are orthogonal to the vectors $\vec
v_{M+k+1},$ $\cdots,$ $\vec v_{M+2k}$ above. This means the 
relations (2.7).

3) The statement can be deduced  from  2).
{\bf Q.E.D.}

\par
In view of the Proposition ~\ref{prop21},  we introduce the subsets
of indices $a \in \{1,2, \cdots, M\}$ as follows.
\begin{dfn}
The subset $I^+_q \subset \{1,2, \cdots, k\}$ (resp. $I^-_q, I^0_q$)
consists of the indices
$a$ such that the coefficient $B^a_q$ of
${\mathcal L}_a(\bi, \bz, \bz', \bzeta)$ (2.5)
is positive (resp. negative, zero).  Analogously we define the subset
$J^+_r \subset \{1,2, \cdots, m\}$ (resp. $J^-_r, J^0_r$)
that consists in such indices
$a$ that the coefficient $C^a_r$  of
${\mathcal L}_a(\bi, \bz, \bz', \bzeta )$
is positive (resp. negative, zero).
\label{dfn2}
\end{dfn}

To assure the convergence of the Mellin inverse transform of 
$M_{{\bi}, \gamma}^\zeta(\bz, \bz')$ from (2.1) 
in a properly chosen angular sector
in the variables $(s,s') \in \bC^{k+m}$,
we shall verify that the Mellin transform $M_{{\bi},\gamma}^\zeta(\bz, \bz')$
admits the following estimation modulo multiplication by a $\Delta-$periodic 
function $g(\bz,\bz').$
$$\mid M_{{\bi},\gamma}^\zeta(\bz, \bz')\mid < C_{\bi}
exp(-\epsilon \mid Im\; z \mid )\;\;
{\rm
while}\;  Im\; z \rightarrow \infty,
\mbox{ in a sector of aperture}\;<2
\pi. $$ for certain $\epsilon >0,$

Here we remember an elementary lemma for the integral:
$$ \int_{z_0 - i\infty}^{  z_0 + i\infty } s^z g(z)
\prod_{j=1}^{\nu}\frac {\Gamma ( z+ \alpha_j)}{\Gamma ( z+ \rho_j)}
dz . \leqno(2.12)$$

\begin{lem}
If one chooses one of the following functions $g^{+}(z)$ (resp. $g^{-}(z)$) 
in terms of  $g(\bz, \bz'),$
then the integrand of (2.12) is exponentially decaying as
$ Im \;z $ tends to  $\infty $ within the sector $0 \leq arg \;z < 2\pi,$
(resp. $-\pi \leq arg \;z < \pi .$)

$$ g^{\pm}(z)= 1+ e ^{\pm 2 \pi i \beta_{\nu}}\prod_{j=1}^{\nu}
\frac {sin 2 \pi ( z+ \alpha_j)}{
sin 2\pi ( z+ \rho_j)}, $$
with $ \beta_{\nu} =-1 + \sum_{j=1}^{\nu}(\rho_j  - \alpha_j)  $
\label{lem221}
\end{lem}

{\bf Proof }

It is enough to recall
$$\prod_{j=1}^{\nu}\frac {\Gamma ( x+ iy + \alpha_j)}{\Gamma ( x+iy + \rho_j)}
\rightarrow const. \mid y \mid^{-(\beta_{\nu}+1)}$$
while $y \rightarrow \pm \infty.$ Here we used the formula of Binet:
$$ log\; \Gamma(z+a) = \log\; \Gamma(z) + a \log\; z - \frac{a- a^2}{2z} + {\cal O}
( \mid z \mid^{-2}   )
$$
if $\mid z \mid >> 1,$ .
The factor $\mid s^{-(x+iy)} \mid = r^{-x}e^{\theta y},$
for $s =  re^{i \theta} $ gives the exponentially decreasing 
contribution in each cases.
{\bf Q.E.D.}

Let us introduce a simplified notation,
$$  {\mathcal L}_j(z)=
A_{j1}z_1 + A_{j2}z_2 +\cdots + A_{jk}z_k + A_{j0}, \; 1 \leq j \leq p,$$
$${\mathcal M}_j (z)=B_{j1}z_1 + B_{j2}z_2 +\cdots + B_{jk}z_k + B_{j0},
\; 1 \leq j \leq r.$$
\begin{lem}
The sufficient conditions so that
$$ \int_{\check \Pi} s^{\bz} g(z)
\frac {\prod_{j=1}^{p}\Gamma ({\mathcal L}_j(z))}{\prod_{j=1}^{r}
\Gamma ( {\mathcal M}_j(z))} dz_1\wedge \cdots \wedge dz_k \leqno(2.13)$$
defines a polynomially increasing function with
 $g(z)$ a properly chosen $\Delta-$periodic function 
(including the infinity $\infty$) are the following.

i) For every $i>0$
$$\sum_{j=1}^p A_{j,i}=\sum_{j=1}^r B_{j,i}  $$

ii)
The real number
$$\alpha = min_{z \in S^{k-1}}\bigl(\sum_{j=1}^p |{\mathcal L}_j(z)- A_{j0}|-
\sum_{j=1}^r|{\mathcal M}_j(z)- B_{j0}| \bigr) $$
is non negative.

\label{lem}
\end{lem}
To see the exponential decay property of the integrand,
one shall make reference to  N\"orlund's trick
\cite{Nor}.
Further we apply the Stirling's formula on the asymptotic behaviour of the $\Gamma-$function (Whittaker-Watson, Chapter XII, Example 44).

If we apply this lemma to our integral, we see that there exists a cycle 
$\check \Pi$  such that
$$I_{x^{\bi}, \gamma}^{\zeta} (s,s'):
= \int_{\check \Pi} g(\bz, \bz')\frac{ \prod_{a\in I^+_q \cup I^0_q}\Gamma
\bigl({\mathcal L}_a (\bi , \bz, \bz' , \zeta)\bigr)} {\prod_{\bar{a}\in
I^-_q }\Gamma \bigl (1-{\mathcal L}_{\bar{a}} (\bi , \bz, \bz' ,
\zeta)\bigr)} s^{-\bz} s'^{-\bz'}d\bz \wedge d\bz', \leqno(2.14)$$ 
with a $\Delta-$periodic function
$g(\bz, \bz')$ rational with respect to  
$e^{2 \pi \sqrt -1 {\mathcal L}_a (\bi , \bz, \bz' ,
\zeta)}, a \in I$. Here we remember the relation
$e^{\pi \sqrt -1 z}\Gamma(z)\Gamma(1-z)= \frac{\pi} {1- e^{-2\;\pi \sqrt -1 z}}.$ 
Thus we get the theorem on the Horn type system.
\begin{thm}
The integral $I_{x^{\bi}, \gamma}^{\bzeta}(s,s')$ satisfies the hypergeometric
system of Horn type as follows:
$$ L_{q ,\bi}(\vartheta_s,\vartheta_{s'} s, s', \zeta)
{I}_{x^{\bi}, \gamma}^{\zeta}(s,s'):=
\Bigl[P_{q, \bi}( \vartheta_s, \vartheta_{s'}, \zeta)
-s_q^\Delta Q_{q, \bi}( \vartheta_s,\vartheta_{s'}, \zeta)\Bigr]
{I}_{x^{\bi}, \gamma}^{\zeta}(s,s')=0,1 \leq q \leq k \leqno(2.15)_1$$
with
$$P_{q, \bi}(\vartheta_s, \vartheta_{s'}, \zeta)=
 \prod_{a\in I^+_q} \prod_{j=0}^{B_q^a-1}
\bigl({\mathcal L}_a(\bi,-\vartheta_{s}, -\vartheta_{s'},
\zeta)+j\bigr), \leqno(2.15)_2$$ 
$$Q_{q, \bi}( \vartheta_s,
\vartheta_{s'}, \zeta) = \prod_{\bar{a}\in
I^-_q}\prod_{j=0}^{-B_q^{\bar a}-1} \bigl({\mathcal
L}_{\bar{a}}(\bi,-\vartheta_{s},-\vartheta_{s'}, \zeta)+j \bigr)
,\leqno(2.15)_3$$
where $I^+_q, I^-_q,  1 \leq q \leq k $
are the sets of indices defined in  Definition ~\ref{dfn2}.
$$ L'_{r ,\bi}(\vartheta_s,\vartheta_{s'}, s, s', \zeta)
{I}_{x^{\bi}, \gamma}^{\zeta}(s,s'):=
\Bigl[P'_{r, \bi}( \vartheta_s, \vartheta_{s'}, \zeta)
-s_r'^\Delta Q'_{r, \bi}( \vartheta_s,\vartheta_{s'}, \zeta)\Bigr]
{I}_{x^{\bi}, \gamma}^{\zeta}(s,s')=0,1 \leq q \leq k \leqno(2.15)_4$$
$$P'_{r, \bi}(\vartheta_s, \vartheta_{s'}, \zeta)=
 \prod_{a\in J^+_r} \prod_{j=0}^{C_r^a-1}
\bigl({\mathcal L}_a(\bi,-\vartheta_{s}, -\vartheta_{s'}, \zeta)+j\bigr)
\leqno(2.15)_5$$
$$Q'_{r, \bi}( \vartheta_s, \vartheta_{s'}, \zeta)
= \prod_{\bar{a}\in J^-_r}\prod_{j=0}^{-C_r^{\bar a}-1}
\bigl({\mathcal L}_{\bar{a}}(\bi,-\vartheta_{s},-\vartheta_{s'},
\zeta)+j \bigr). \leqno(2.15)_6$$
where $J^+_r, J^-_r,  1 \leq r \leq m$ are the sets of indices defined in the
Definition ~\ref{dfn2}.
The degree of two operators  $P_{q,
\bi}(\vartheta_s, \vartheta_{s'}, \zeta),$ $Q_{q, \bi}( \vartheta_s,\vartheta_{s'}, \zeta)$ are equal.
Namely,
$$ deg\;P_{q,
\bi}(\vartheta_s, \vartheta_{s'}, \zeta) =  \sum_{a \in I^+_q} B_q^a = - \sum_{\bar a \in I^-_q} B_q^{\bar a} =deg\;Q_{q,
\bi}(\vartheta_s, \vartheta_{s'}, \zeta) .
\leqno(2.16)$$ Analogously,
$$ deg\;P'_{r,
\bi}(\vartheta_s, \vartheta_{s'}, \zeta) =  
\sum_{a \in J^+_r} C_r^a = - \sum_{\bar a \in J^-_r} C_r^{\bar a} =
deg\;Q'_{r, \bi}( \vartheta_s, \vartheta_{s'}, \zeta).$$

\label{thm23}
\end{thm}

The proof is mainly based on the Proposition ~\ref{prop21}.
To deduce $(2.15)$ from the Mellin transform 
$M_{\bi, \gamma}^\zeta(\bz, \bz')$
we use the
following well known recurrence relation:
$$\Gamma(\frac{\alpha(n+\Delta)}{\Delta}+\zeta)=
\Gamma(\frac{\alpha n}{\Delta}+\zeta)(\frac{\alpha n}{\Delta}+\zeta)
(\frac{\alpha n}{\Delta}+1+ \zeta)\cdots (\frac{\alpha n}{\Delta}+ \alpha-1+\zeta),
$$ if $\alpha>0$ a positive integer.
$$\Gamma(\frac{\alpha(n+\Delta)}{\Delta}+\zeta)=
\Gamma(\frac{\alpha n}{\Delta}+\zeta)(\frac{\alpha
n}{\Delta}+\zeta-1)^{-1} (\frac{\alpha n}{\Delta}+
\zeta-2)^{-1}\cdots (\frac{\alpha n}{\Delta}+ \zeta+ \alpha)^{-1},
$$ if $\alpha<0$ a negative integer.

The evident compatibility (i.e. integrability) of the above system
$(2.15)_\ast$ in the sense of Ore-Sato
(\cite{Sad1}) can be formulated like the following cocycle condition.
To state the proposition
we introduce the notation $\bz +\Delta e_r = (z_1, \cdots, z_{r-1}, 
z_r+\Delta, z_{r+1}, \cdots, z_k).$
\begin{prop}
The rational expression
$$R_q(\bz, \bz')= \frac{P_{q, \bi}(\bz,\bz', \zeta)}{Q_{q, \bi}
(\bz+\Delta e_q, \bz',\zeta )}, \leqno(2.17)$$ defined for the operators
$(2.15)_2,$ $(2.15)_3$ satisfies the following relation:
$$R_q(\bz+\Delta e_r,\bz')R_r(\bz, \bz')= R_r(\bz+\Delta e_q, \bz')
R_q(\bz, \bz'), \;\; q,r =1, \cdots, k.
\leqno(2.18) $$
Similarly  for
$$R'_\kappa(\bz, \bz')= 
\frac{P'_{\kappa, \bi}(\bz,\bz', \zeta)}{Q'_{\kappa, \bi}
(\bz, \bz'+ \Delta e'_\kappa,
\zeta )}, \leqno(2.19)$$ satisfies the following relation:
$$R'_\kappa(\bz, \bz'+\Delta e'_\rho)R'_\rho(\bz, \bz')= 
R'_\rho(\bz, \bz'+\Delta e'_\kappa)
R'_\kappa(\bz, \bz'), \;\; \kappa,\rho =1, \cdots, m.
\leqno(2.20) $$
\end{prop}

\begin{remark}
As $m = dim \;{\bf D}(\tilde \Sigma)$(see (1.3)), one can consider that the above system $(2.15)_{\ast}$ is defined on 
 $ \bT^k \times {\bf D}(\tilde \Sigma)$
for ${\bf D}(\tilde \Sigma)$ : the 
N\'eron-Severi torus associated to the fan 
$\tilde \Sigma$. \label{remark411}
\end{remark}
We introduce here the main object of our study: the
discriminantal loci of the CI defined by the polynomials
$f_1(x,{s}') +s_1,\cdots, f_k(x,{s}') +s_k.$
$$ D_{s,s'} := \{(s,s') \in \bT^{k+m}; f_1(x,{s}') +s_1=\cdots=
f_k(x,{s}') +s_k=0,
rank
\left (
\begin {array}{c} grad_x \; f_1(x,{s}')\\
\cdots \\ grad_x \; f_k(x,{s}')\\
\end{array}
\right) < k, \leqno(2.21) $$
\begin{flushleft}
$\rm{for \; certain}\;$ $x\in
\bT^N\}.$
\end{flushleft}
As it is easy to see \cite{GZK2}, $D_{s,s'}$ coincides with the discriminantal loci
of $F(x,s',s,y).$

Let us define the $\Delta-$th roots of rational functions associated with the linear functions $(2.5)$ as follows.
$$ \psi_q(\bz, \bz') = \left(\frac{\prod _{a \in I_q^+}(\sum_{\ell=1}^k
B_\ell^a z_\ell
+\sum_{j=1}^{m} C_j^a z'_j)^{B_q^a}} {\prod _{\bar a \in I_q^-}
(\sum_{\ell=1}^k  B_\ell^{\bar a} z_\ell+ \sum_{j=1}^{m} C_j^{\bar
a} z'_j)^{-B_q^{\bar a}}}\right)^{\frac{1}{\Delta}}, 1 \leq q \leq
k, \leqno(2.22)$$
$$ \phi_r(\bz, \bz') = \left(\frac
{\prod _{a \in J_r^+}(
\sum_{\ell=1}^k B_\ell^a z_\ell
+\sum_{j=1}^{m} C_j^a z'_j)^{C_r^a}}
{\prod _{\bar a \in J_r^-}(\sum_{\ell=1}^k  B_\ell^{\bar a} z_\ell
+\sum_{j=1}^{m} C_j^{\bar a} z'_j )^{-C_r^{\bar a}}
}\right)^{\frac{1}{\Delta}}, 1 \leq r \leq m. \leqno(2.23)
$$ 
$$
h: \bC^{k+m}\setminus \{0\} \rightarrow (\bC^{\times})^{k+m}, \leqno(2.24)$$
$$(\bz,\bz') \rightarrow (\psi_1(\bz,\bz'), \cdots ,
\psi_k(\bz,\bz'), \phi_1(\bz,\bz'), \cdots ,\phi_m(\bz,\bz')).$$
By virtue of the property  $(2.7),$ the rational function
$\psi_q(\bz,\bz')^\Delta$
(resp. $\phi_r(\bz,\bz')^\Delta$) is of weight zero with respect 
to the variables 
$(\bz,\bz')$ and thus it is possible to consider the mapping $h$
defined on $\bC P^{k+m-1}$ instead of $\bC ^{k+m}.$

Let $\Delta_f(s,s')$ be a polynomial that defines the discriminantal loci
$D_{s,s'}$ without multiplicity.
\begin{thm} The image of $
h: {\bC \bP}^{k+m-1}\rightarrow (\bC^{\times})^{k+m}$
is identified with the discriminantal loci  $D_{s,s'}$
if we choose a proper $\Delta-$th branch in the equations $(2.2),$
$(2.3).$  
\label{thm22}
\end{thm}
{\bf Proof}
From the system of equations $(2.15)$ we see that 
$D_{s,s'}$ is contained in the set:
$$ \nabla_{s,s'} := \{(s,s') \in \bT^{k+m};  \sigma(L_{q,-\b1})
(s\xi, s'\xi',s,s', -\b1)
=0,1 \leq q \leq k, \leqno (2.25) $$
$$ \sigma(L'_{r,-\b1})(s\xi, s'\xi',s,s', -\b1)
=0,1 \leq r \leq m \;\;{\rm for\; some \;} 
(\xi, \xi') \in \bT^{k+m} \}. $$
here we use the notation
$(s\xi, s'\xi')$ $= (s_1 \xi_1, \cdots, s_k \xi_k, s'_1 \xi'_1, 
\cdots, s'_m \xi'_m ).$
The existence of $(\xi, \xi') \in \bT^{k+m}$ in $(2.25)$
is equivalent to the existence of $(\bz,\bz')= (s\xi, s'\xi') 
\in \bT^{k+m}.$ Thus the set $ \nabla_{s,s'} $ admits a representation,
$$ \{(s,s') \in \bT^{k+m}; s_q^\Delta =\frac{P_{q,-\b1}
(\bz, \bz', -\b1)}{Q_{q,-\b1}
(\bz, \bz', -\b1)} \ , 1 \leq q \leq k, \;  (s'_r)^\Delta
= \frac{P'_{r,-\b1}
(\bz, \bz', -\b1)}{ Q'_{q,-\b1}
(\bz, \bz', -\b1)}, 1 \leq r \leq m \}. $$
While after Theorem 2.1, a) and Remark 2.4 of \cite{Kap1}, this set
$\nabla_{s,s'} $ coincides with $D_{s,s'}$ if $\Delta=1.$
As for the case $\Delta >1,$ it is natural to consider the $\Delta-$covering
$\tilde h$ of the mapping $h,$
$$ \tilde h: {\bC \bP}^{k+m-1}\rightarrow \tilde{(\bC^{\times})^{k+m}},$$
while the branch of the image of $h$ shall be specified in a proper way.
To do that we remark that $ h({\bC \bP}^{k+m-1})$ $\subset \nabla_{s,s'}$
where the difference $ \nabla_{s,s'} \setminus h({\bC \bP}^{k+m-1}) $ 
consists of the divisors that arise from the $\Delta-$branching effect $ \tilde h({\bC \bP}^{k+m-1}).$ In considering $D_{s,s'}$ we shall discard the 
superfluous $\Delta-$branching effect $ \tilde h({\bC \bP}^{k+m-1}) \setminus h({\bC \bP}^{k+m-1}).$ 
{\bf Q.E.D.}

The mapping $(2.24)$ is nothing but the inverse mapping of the logarithmic
Gauss map;
$$ D_{s,s'}
\rightarrow \bC P^{k+m-1},$$
$$ (s,s') \rightarrow \bigl(
s_1\frac{\partial}{\partial s_1}\Delta_f(s,s'): \cdots :
s_k\frac{\partial}{\partial s_k}\Delta_f(s,s'):
s_1'\frac{\partial}{\partial s_1'}\Delta_f(s,s'): \cdots:
s_m'\frac{\partial}{\partial s_m'}\Delta_f(s,s') \bigr).$$
This is a direct consequence of the cocycle property $ (2.18), $
$(2.20)$
of the operators $L_{q ,\bi}(\vartheta_s,\vartheta_{s'}, s, s', \zeta)$
and $L'_{r ,\bi}(\vartheta_s,\vartheta_{s'}, s, s', \zeta),$ see 
\cite{Kap1}, Theorem 2.1, b).

\vspace{2pc}
{
\center{\section{$A-$Hypergeometric function of Gel'fand-Kapranov-Zelevinski
 }}
}

Let us consider the set of polynomials with deformation parameter coefficients
$(a_{0,1}, \cdots, a_{\tau_k,k}) $ associated to the polynomial system
$(0.2),$ 
$${\bar f}_\ell(x, \ba)=
a_{1,\ell}x^{\vec \alpha_{1,\ell}}+\cdots +a_{\tau_\ell,\ell}
x^{\vec \alpha_{\tau_\ell,\ell}}+a_{0,\ell}. \; 1 \leq \ell \leq k.
\leqno(3.1) $$
For the sake of simplicity we will further make use of the notation
$\ba:=(a_{0,1}, \cdots, a_{\tau_k,k}) \in \bT^L.$
We consider the Leray coboundary $\partial \gamma_\ba$ of a cycle
$\gamma_\ba \in H_{n}(X_{\ba},\bZ)$ of the CI $X_\ba=\{x \in \bT^N;
{\bar f}_1(x, \ba)= \cdots ={\bar f}_k(x, \ba)=0 \}.$

Then we can define the $A-$ hypergeometric function $\Phi_{x^\bi, \gamma_{\ba}}^\zeta(a_{0,1}, \cdots, a_{\tau_k,k})$ introduced by
  Gel'fand-Zelevinski-Kapranov
\cite {GZK} associated to the polynomials,
$$ f_\ell(x) = x^{\vec \alpha_{1,\ell}}+\cdots + x^{\vec \alpha_{\tau_\ell,\ell}},
\; 1 \leq \ell \leq k,$$
$$ x^{\bi}= x_1^{i_1}\cdots x_N^{i_N}, x^{\vec \alpha_{j,\ell}}=
x_1^{\alpha_{j,\ell,1}} \cdots x_N^{\alpha_{j,\ell,N}}.$$
Namely it is defined as a kind of multiple residue along $X_\ba,$
$$\Phi_{x^\bi, \gamma_{\ba}}^\zeta(a_{0,1}, \cdots, a_{\tau_k,k}):=
\int_{\partial \gamma_{\ba}} \prod_{\ell=1}^k
 {\bar f}_\ell(x, \ba)^{-\zeta_\ell-1}
x^{\bi +\b1} \frac{dx}{x^{\b1}}. \leqno(3.2)$$
We impose here the non-degeneracy condition  
of the Definition 
~\ref{dfn15} for the complete intersection 
${X}_s$ after the procedure described in \S 1.

In the sequel we consider a lattice $\Lambda \subset \bZ^L$
of $L-$vectors defined by the system of following linear equations:
$$ \sum^{\tau_q}_{i=0}b(j,q,\nu)=0,\;\; 1 \leq q \leq
k,$$ $$\sum_{q=1}^k\sum_{j=1}^{\tau_q}\alpha_{j q \ell}
b(j,q,\nu)=0,  1 \leq \ell \leq N.$$
Here we denoted by $(b(0,1,\nu), \cdots,b(\tau_1,1,\nu),
b(0,2,\nu),\cdots, b(\tau_2,2,\nu),\cdots, b(\tau_k,k,\nu)),$
$1 \leq \nu \leq
m+k,$ a $\bZ$ basis of $\Lambda.$

For the subset  $\bK \subset \{(0,1),\cdots, (k,\tau_k)\}$
such that the columns  ${\vec m_{j,q}}(A), (j,q) \in \bK$ of the matrix
${\sf M}(A)$ (1.7) span $\bR^{N+k}$ over  $\bR$ and $|\bK|= N+k$
we define the set of indices (a generalisation of the Frobenius' method)
after \cite{GZK},
$$ \Pi((\zeta+\b1, \bi +\b1), \bK)=\{((\lambda(0,1,\nu), \cdots,
\lambda(\tau_1,1,\nu),
\cdots, \lambda(\tau_k,k,\nu))\}_{
1 \leq \nu \leq |det ( {\vec m_{j,q}}(A))_{(j,q) \in \bK}| },$$
which satisfy the following system of equations,
$$ \sum^{\tau_\nu}
_{j=0}\lambda(j,q,\nu) +\zeta_q +1
=0,\;\; 1 \leq q \leq k,$$
$$\sum_{q=1}^k\sum_{j=1}^{\tau_q}
\alpha_{j q \ell}\lambda(j, q, \nu)
-(i_\ell+1)=0,  1 \leq \ell \leq N.$$  Let $T$ be a triangulation of the 
Newton polyhedron $\Delta(F(x,\b1, \b1,y)+1)$ for  $F(x,\b1, \b1, y)$ of (1.4)
after the definition \cite{GZK}, 1.2. Here we impose that
 $ \lambda(j,q,\nu) \in \bZ $ for $ (j,q) \not \in \bK.$
Let $\bK_1, \bK_2 \in T$ be two different simplices
of the 
triangulation $T.$ We suppose that $\vec \lambda(\nu_p):= (\lambda(0,1,\nu_p), \cdots,
\lambda(k,\tau_k, \nu_p)) \in \Pi((\zeta+\b1, \bi +\b1 ), \bK_p),$
$\lambda(j,q,\nu_p) \in
\bZ$ for $(j,q) \not \in \bK_p,$ $(p=1,2)$ with $ 1 \leq \nu_p \leq |det (
  {\vec m_\rho}(A))_{\rho \in \bK_p}|.$ We introduce the condition of 
$T-$ non-resonance on $(\zeta+\b1, \bi+\b1)$ $$
(\lambda(0,1,\nu_1), \cdots,
\lambda(k,\tau_k, \nu_1))\not \equiv
(\lambda(0,1,\nu_2), \cdots,
\lambda(k,\tau_k, \nu_2))
\; mod \;\Lambda,  \leqno(3.3)$$ for any pair
$ \vec \lambda(\nu_p)=(\lambda(0,1,\nu_p), \cdots,
\lambda(k,\tau_k, \nu_p))\in \Pi((\zeta+\b1, \bi +\b1 ),
\bK_p),$ $p=1,2.$ An adaptation of theorem 3 \cite{GZK} \
to our situation can be formulated as follows.

\begin{thm}
1) The $A-$HGF $\Phi_{x^\bi, \gamma_a}^\zeta(\ba)$
satisfies the following system of 
equations.

$$\bigl(\sum_{j=0}^{\tau_q}a_{{ji}}\frac{\partial}{\partial a_{ji}}
+\zeta_q +1 \bigr)\Phi_{x^\bi, \gamma_a}^\zeta(\ba)=0, 1 \leq q \leq k,
\leqno(3.4) $$
$$
\bigl(\sum_{1 \leq q \leq k, 1 \leq j \leq \tau_q}
\alpha_{jq1}a_{{jq}}\frac{\partial}{\partial a_{jq}}
-(i_1+1) \bigr)\Phi_{x^\bi, \gamma_a}^\zeta(\ba)=
\cdots$$ $$ =\bigl(\sum_{1 \leq q \leq k, 1 \leq j \leq \tau_q}
\alpha_{jqN}a_{{jq}}\frac{\partial}{\partial a_{jq}}
-(i_N+1) \bigr)\Phi_{x^\bi, \gamma_a}^\zeta(\ba)= 0,$$
$$ \bigl(\prod_{\{(j,q);b(j,q,\nu)>0\}}
(\frac{\partial}{\partial a_{jq}})^{b(j,q,\nu)}
- \prod_{\{(j,q);b(j,q,\nu)<0\}}
(\frac{\partial}{\partial a_{jq}})^{-b(j,q,\nu)}\bigr)
\Phi_{x^\bi, \gamma_a}^\zeta(\ba)= 0,\; 1 \leq \nu \leq L-(k+N).$$

2) The dimension of solutions of the system above
at a generic point  $\ba \in
\bT^L$ is equal to $(N+k)! vol_{N+k}
\Delta (F(x,\b1, \b1, y)+1)$ $= |\chi ( Z_{F(x, \b1,\b1, y )} )|$ 
if the $T-$non-resonant 
condition $(3.3)$ is satisfied. \label{thm411}
\end{thm}

In the sequel we  shuffle  the variables $\ba
=(a_{0,1}, \cdots, a_{\tau_k,k})$ in accordance with  the order of 
their appearance
and we define anew the indexed parameters
$a_{1}=a_{1,1}, \cdots,$ $a_{\tau_1}=a_{\tau_1,1},$ $a_{\tau_1+1}=
a_{0,1}, \cdots,$ $a_{L-1}=a_{\tau_k,k},$ $a_{L}=a_{0,k}.$
Let us introduce notations analogous to  $(1.14),$
$$\Xi (A):= ^t(\log\; X_1, \cdots, \log\; X_N,  \log\; a_1, \cdots,
log\; a_L,\log\; U_1, \cdots,
log\; U_k). \leqno(3.5)$$
$$ \log\; T_1 = <\vec \alpha_{1,1}, \log\;X> + log\; a_1 + \log\; U_1,$$
$$\vdots$$
$$ \log\; T_{\tau_1} = <\vec \alpha_{1,\tau_1}, \log\;X> + log\; a_{\tau_1}
+ \log\; U_1,$$
$$\vdots$$
$$ \log\; T_L = log\; a_L  + \log\; U_k.$$

We consider the equation
$$ {\sf L}(A)\cdot Log\; \Xi(A) = {\sf L}\cdot Log\; \Xi,$$
where the matrix  ${\sf L}(A)$ is constructed as follows.
The columns  $\vec \ell_i(A) = \vec v_i, 1 \leq i \leq N$
with vectors $\vec v_i$ defined like the column of the matrix
$\sf L$ in (1.15). For the columns of number $N+1$ to $N+L$
$(\vec \ell_{N+1}(A), \cdots, \vec \ell_{N+L}(A))= id_L.$
The columns $\vec \ell_{N+L+j}(A)=
^t(\overbrace{0,\cdots, 0,0,}^{\tau_1+\cdots +\tau_{j-1}+j-1}
,\overbrace{1,1, \cdots, 1}^
{\tau_{j}+1}, 0, \cdots,0), 1 \leq j \leq k,$
 the matrix ${\sf L}(A)$ is obtained after implementation of the matrix $id_L$ into the transposed matrix $^t{\sf M}(A)$
between the $k-$th and the $(k+1)-$th column up to necessary permutations necessary after the implementation.

 \begin{prop} There exists a cycle $\gamma_a$ such that the following equality
holds for the integral defined in
$(3.2),$ $$ \Phi_{x^\bi, \gamma_\ba}^\zeta(\ba) = B_\bi^\zeta(\ba)
I_{x^\bi, \gamma}^\zeta(s(\ba), s'(\ba)), \leqno(3.6)$$ here $$
s_\ell(\ba)= \prod_{j=1}^L a_j^{w_{j, N+\ell}}, \;\;, 1\leq \ell
\leq k,$$ $$ s_\rho'(\ba)= \prod_{j=1}^L a_j^{w_{j, N+k+\rho}},
\;\;, 1\leq \rho \leq m,$$ $$B_\bi^\zeta(\ba)= \prod_{\ell=1}^N
(\prod_{j=1}^L a_j^{w_{j,\ell}})^{i_\ell+1 }
\prod_{\nu=1}^k (\prod_{j=1}^L a_j^{w_{j,N+k+m+\nu}})^{\zeta_\nu+1}.$$
The exponents $w_{j,\ell}$ 
are determined by the following relation,
$${\sf L}^{-1}\cdot {\sf L}(A)
= \left [\begin {array}{ccccccccc}
1&\cdots&0&w_{1,1}&\cdots&w_{L,1}&0&&\\
\vdots&\ddots&\vdots&\vdots&\vdots&\vdots&\vdots&\vdots&\vdots\\
0&\cdots&1&w_{1,N}& \cdots& w_{L,N}&0&&\\
0&\cdots&0&w_{1,N+1}& \cdots& w_{L,N+1}&0&&\\
\vdots&\ddots&\vdots&\vdots&\vdots&\vdots&\vdots&\vdots&\vdots\\
0&\cdots&0&w_{1,N+k+m}& \cdots& w_{L,N+k+m}&0&&\\
0&\cdots&0&w_{1,N+k+m+1}& \cdots& w_{L,N+k+m+1}&1&\cdots&0\\
\vdots&\ddots&\vdots&\vdots&\vdots&\vdots&\vdots&\vdots&\vdots\\
0&\cdots&0&w_{1,L}& \cdots& w_{L,L}&0&\cdots&1\\
\end {array}\right ],
\leqno(3.7)$$
that has been essentially introduced in $(2.8)$. 
The transition of the cycle $\gamma(a)$ to
$\gamma$ is controlled by the transformations,
$$ X_i = (\prod_{j=1}^L a_j^{w_{j,i}})^{-1}
 \cdot x_i.$$
\label{prop42}
\end{prop}
{\bf Proof}
It is enough to remark the following property,
$$ x^{\bi+\b1} y^{\zeta +\b1} \frac{dx}{x^\b1}\wedge \frac{dy}{y^\b1}
=  B_\bi^\zeta(\ba) x^{\bi+\b1} U^{\zeta +\b1} \frac{dx}{x^\b1}\wedge
\frac{dU}{U^\b1}.$$
{\bf Q.E.D.}
\par One can thus conclude (at least locally on the chart
$a_j \not = 0$ for $j \in I, |I|= k+m$)
$A-$HGF of GZK $(3.2)$ is expressed by means of a fibre integral
annihilated by the Horn system $(2.15).$
One can find a similar statement in \cite{Kap1}
where Kapranov restricts himself to a power series expansion of
the solution to $(3.2)$.

\begin{cor}
The dimension of the solution space of the system
(3.3) at the generic point is equal to 
$ |\chi ( Z_{F(x, \b1, \b1, y )} )|$ if the $T-$
 non-resonance condition $(3.3)$ is satisfied. \end{cor}

{\bf Proof} 
We shall consider the convex hull of vectors that correspond to the vertices of the Newton polyhedron of the polynomial $ y_1(f_1(x)+1)+$
$\cdots$ $+y_k(f_k(x)+1).$ That is to say $(\vec \alpha_{1,1},1,$ $ 0,$
$\cdots,0),$ $\cdots,$ $(\vec \alpha_{\tau_1,1},$ $1, 0,\cdots,0),$
$(\vec \alpha_{1,2},0,1, 0, \cdots,0),$ $\cdots,$ $(\vec \alpha_{\tau_k,k},0,
\cdots,0,1) \in \bZ^{N+k}.$ They are located on the hyperplane $\zeta_1 +
\cdots +\zeta_k=1.$
Thus it is possible to measure
$(N+k-1)$ dimensional volume $(N+k-1)! vol_{N+k-1}(\Delta(F(x,\b1, \b1, y))$ that is equal to $(N+k)! vol_{N+k}(\Delta(F(x,\b1, \b1, y)+1).$
The Euler characteristic admits the following expression
$$ |\chi(Z_{F(x,\b1,\b1,y)})|= \sum_p |det{\sf M}_{\bK_p}|=(N+k-1)!
vol_{N+k-1}\bigl(\Delta(F(x,\b1,\b1,y)
)\bigr),$$
after Khovanski \cite{Kh1}. {\bf Q.E.D.}

We define the $A-$discriminantal loci $ \nabla_{\ba}^0$ in $\bT^L$
like following,
$$ \nabla_{\ba}^0 = \{\ba \in \bT^{L}; \bar f_1(x,\ba) =\cdots=\bar
f_k(x,\ba) =0,
rank
\left (
\begin {array}{c} grad_x \; \bar f_1(x,\ba)\\
\cdots \\ grad_x \; \bar f_k(x,\ba)\\
\end{array}
\right) < k \}. \leqno(3.8)$$
As it is seen from  $(3.7)$ the uniformisation equations $(2.22),$ $(2.23)$
give rise to an uniformisation of $A-$discriminantal loci $ \nabla^0_{\ba}$
without $\Delta-$ branching effect. 

\begin{cor}
We have the following relations among $\ba \in \bT^L$ located on the discriminantal loci $ \nabla_{\ba}^0,$
$$ \prod_{j=1}^L \bigl(\frac{a_j}{
{\mathcal L}_j(-\b1,\bz, \bz',-\b1)}\bigr)^{B^q_j}=1,  \;\; 1 \leq q \leq k,
\leqno(3.9)_1$$
$$ \prod_{j=1}^L \bigl(\frac{a_j}{
{\mathcal L}_j(-\b1,\bz, \bz',-\b1)}\bigr)^{C^r_j}=1,  \;\; 1 \leq r \leq m.
\leqno(3.9)_2$$
This allows us to express $ \nabla_{\ba}^0$ by means of the deformation 
parameters $(\bz, \bz')$ $\in \bC\bP^{k+m-1}$ and  $\ba' $
$\in {\bT^{L-k}}/{\bD(\Sigma)}$ $\cong \bT^{L-(k+m)}.$
\end{cor}
\vspace{2pc}
{
\center{\section{Examples
 }}
}

{\bf 4.1} Deformation of $D_4.$ 

Let us consider the versal deformation of $D_4$ singularity of the
following form,
$$ f(x,s_0,s_1,s_2,s_3)=x_1^3 + x_1 x_2^2 + s_3 x_1^2 + s_2 x_1 + s_1 x_2 + s_0. \leqno(4.1) $$
By means of the resultant calculus on computer, we get a defining equation of the discriminantal loci as follows,
$$ \Delta_f(s)= 1024 s_1^6(432 s_0^4  + 64 s_1^6  + 576 s_0^2  s_1^2  s_2 + 128 s_1^4  s_2^2  + 64 s_0^2  s_2^3  + 64 s_1^2  s_2^4  + 192 s_0 s_1^4  s_3  \leqno(4.2)$$
$$- 288 s_0^3  s_2 s_3 - 320 s_0 s_1^2  s_2^2  s_3 - 24 s_0^2  s_1^2  s_3^2  - 144 s_1^4  s_2 s_3^2  -  16 s_0^2  s_2^2  s_3^2$$ $$  - 16 s_1^2  s_2^3  s_3^2  + 64 s_0^3  s_3^3  + 72 s_0 s_1^2  s_2 s_3^3  + 27 s_1^4s_3^4). $$
This is a polynomial with quasihomogeneous weight $24$ if we assign to the variables $(x_1,x_2; s_0, s_1, s_2,s_3)$ the weights $(1,1; 3 ,2 ,2 , 1).$
Here we remark that $s_1=0$ branch  of the discriminantal locus $D_s =
\{s \in \bC^3; \Delta_f(s)=0 \}$ corresponds to the deformation of $A_2$
singularity.
  
On the other hand, our Theorem 2.6 states that the uniformisation equation of 
the discriminantal loci for the deformation (i.e. torus action quotient of 
the deformation parameter space $(s_0,s_1,0,s_3)$ on the chart 
$s_3 \not = 0$), 
$$ f(x,s_0,s_1,0,1) =x_1^3 + x_1 x_2^2 + x_1^2 + s_1 x_2 + s_0,$$
has the following form,
$$s_0= -\frac {z_2(3 z_1 + 4 z_2 )^2 }{ 4 (2z_1 + 3z_2 )^3}, \leqno(4.3) $$
$$s_1=(- \frac{z_1 (3 z_1 + 4 z_2 )^3 } {4 (2z_1 + 3z_2 )^4})^{\frac 1 2}. $$
If we eliminate the variables $(z_1, z_2)$ from the expressions
$(4.3),$ we get an equation 
$$ 64 s_0^3  +432 s_0^4  -24 s_0^2  s_1^2  + 27 s_1^4  +192 s_0 s_1^4  + 
64 s_1^6 =0.$$ 
We recall here that our method requires that the expression $yf(x,s)$ contains so much terms as the variables in it. The reason why the value $(s_2,s_3) =(0,1)$has been chosen is of purely technical character. 
In substituting the special value  $(0,1)$ for $(s_2,s_3)$ in $(4.2)$
we get, 
$$ \frac{\Delta_f(s_0,s_1,0,1)}{1024 s_1^6} =  64 s_0^3  +432 s_0^4  -24 s_0^2  s_1^2  + 27 s_1^4  +192 s_0 s_1^4  + 64 s_1^6.$$ 

{\bf 4.2} Deformation of a non-quasihomogeneous complete intersection.

Let us consider the following pair of polynomials
that define a non-degenerate complete intersection $X_s $in
$\bC^2$,
$$f_1= x_1^3 + x_2^2 + s_1, f_2 = x_1^2 + x_2^3 + s_2. \leqno(4.4)$$
The discriminant of this CI in $\bC^2$ can be calculated as follows,
$$ (s_1^3+s_2^2)^3(s_2^3+s_1^2)^3 (800000 + 387420489 s_1^5 - 43740000 s_1 s_2 +$$
$$+ 438438825 s_1^2 s_2^2 + 387420489 s_1^3 s_2^3 + 387420489 s_2^5). \leqno(4.5)$$
Evidently the fibres corresponding to the parameter values on the divisor 
$ (s_1^3+s_2^2)^3(s_2^3+s_1^2)^3=0$ are contained in $\{(x_1,x_2) \in \bC^2;x_1x_2=0\}.$
Thus the discriminant of CI $X_s \cap \bT^2$ is given  by the third factor
of $(4.5).$ After Theorem 2.6, we can find an uniformisation equation
of the discriminantal loci  $D_s$,
$$ s_1 = -(\frac {(4z_1+6z_2)^4(5z_1)^5(6 z_1 + 4 z_2 )^6 }{(9z_1 + 6z_2 )^9 (6 z_1 + 9 z_2 )^6})^{1/5}, \leqno(4.6)$$
$$ s_2 = -(\frac {(4z_1+6z_2)^6(5z_2)^5(6 z_1 + 4 z_2 )^4 }{(9z_1 + 6z_2 )^6 (6 z_1 + 9 z_2 )^9})^{1/5}.$$
If we eliminate the variables $(z_1, z_2)$ from the expressions
$(4.6),$ we get an equation of $\nabla_s,$
$$ (800000 + 387420489 s_1^5 - 43740000 s_1 s_2 
+ 438438825 s_1^2 s_2^2 + 387420489 s_1^3 s_2^3 + 387420489 s_2^5) 
R(z_1,z_2),$$
where $R(z_1,z_2)$ is a polynomial whose Newton polyhedron is contained in a
four sided rectilinear figure with vertices $(0,0),$ $(20,0),$ $(12,12),$
$(0,20).$ This factor contains the image of $\tilde h ({\bf CP}^1)$
outside of $D_s.$


\vspace{\fill}

%

\noindent

\begin{flushleft}
\begin{minipage}[t]{6.2cm}
  \begin{center}
{\footnotesize Indepent University of Moscow\\
Bol'shoj Vlasievskij pereulok 11,\\
 Moscow, 121002,\\
Russia\\
{\it E-mails}:  tanabe@mccme.ru, tanabesusumu@hotmail.com \\}
\end{center}
\end{minipage}
\end{flushleft}

\end{document}